\documentclass[12pt]{amsart}

\usepackage{amssymb,amsmath,amsthm,amscd}
\input{amssym.def}

\newcommand{\eea}{\end{eqnarray}}
\newcommand{\ul}{\underline}

\newtheorem{Theorem}{\sc Theorem}
\newtheorem{Lemma}{\sc Lemma}
\newtheorem{Proposition}{\sc Proposition}
\newtheorem{Corollary}{\sc Corollary}
\newtheorem{Definition}{\sc Definition}

\newtheorem{Conjecture}{\sc Conjecture}

\newcommand{\be}{\begin{equation}}
\newcommand{\ee}{\end{equation}}

\def\ul{\underline}%

\def\C{\mbox{${\mathbb C}$}}
\def\D{\mbox{${\mathbb D}$}}
\def\T{\mbox{${\mathbb T}$}}

\def\B{\mbox{${\mathbb B}$}}
\begin{document}

\title{QUASI-FREE RESOLUTIONS OF HILBERT MODULES}

\begin{abstract}

The notion of a quasi-free Hilbert module over a function algebra
$\mathcal{A}$ consisting of holomorphic functions on a bounded domain
$\Omega$ in complex $m$ space is introduced.  It is shown that 
quasi-free Hilbert modules correspond to the completion of 
the direct sum of a certain
number of copies of the algebra $\mathcal{A}$.  A Hilbert module is said to be
weakly regular (respectively, regular) if there exists a module map
from a quasi-free module with dense range (respectively, onto).  A
Hilbert module $\mathcal{M}$ is said to be compactly supported if there exists a
constant $\beta$ satisfying $\|\varphi f \| \leq \beta \|\varphi \|_X \|f\|$ for
some compact subset $X$ of $\Omega$ and 
$\varphi$ in $\mathcal{A}$, $f$ in $\mathcal{M}$.  It is shown that
if a Hilbert module is compactly supported then it is weakly regular.
The paper identifies several other classes of Hilbert modules which
are weakly regular.  In addition, this result is extended to yield
topologically exact resolutions of such modules by quasi-free ones.

\end{abstract}

\author{\tt Ronald G. Douglas} 
\address{Texas A\&M University, College Station, Texas 77843}
\email{rgd@tamu.edu}
\author{\tt Gadadhar Misra}
\address{Indian Statistical Institute, RV College Post, Bangalore 560 059}
\email{gm@isibang.ac.in}
\thanks{The research of both the authors was supported in part by DST -
NSF S\&T Cooperation Programme. This paper was completed while 
the first named author was visiting the Fields Institute on a 
development leave from Texas A \& M University.}
\date{}
\maketitle
%\noindent

\section{Introduction} 

While the study of linear operators on Hilbert space goes back more
than a hundred years, attempts at understanding more than one operator
are of more recent origin.   Rings of operators were
investigated in the celebrated series of papers by Murray and von
Neumann \cite{Mvon0}, \cite{Mvon1}, \cite{Mvon2}, \cite{von2}
in the thirties, but that study led to the development of
operator algebras. This subject is somewhat different than operator theory,
and most recently has led to noncommutative geometry \cite{AC}.  Also, 
there have been several approaches to non-selfadjoint
multivariate operator theory. For example, there is the study of
non-selfadjoint operator algebras which was initiated by Kadison and
Singer \cite{KS} and has been developed by many authors over the years
(cf.  \cite{KD}).  In \cite{WBA2}, \cite{WBA3}, Arveson extended 
results on function algebras, especially the disk algebra, to 
non-selfadjoint subalgebras of $C^*$-algebras.
There is the study of operators in various concrete 
settings, usually defined on spaces of holomorphic functions such as
the Hardy and Bergman spaces (cf. \cite{AMS}, \cite{NKN}).
Generalizations of operators on Hardy space have been undertaken by 
several authors recently including Arveson \cite{WBA4}, 
Popescu \cite{GP2}, and Davidson \cite{KD1}.
Finally, there is the module approach (cf. \cite{rgdvip}), in
which an algebraic point of view is emphasized and the extension of
techniques from algebraic and complex geometry  is the
key.  This note makes a contribution to the latter program.

In commutative algebra a principal object of study is the collection
of modules over a given ring.  While in most instances the
collection has additional structure and forms a semigroup, a group or
even a ring, one way to study individual modules is by relating them
to simpler ones.  If the ring is Noetherian, then one uses projective
modules which can be characterized in this context as submodules of
free modules.  The latter are the direct sum of copies of the
ring.  Otherwise, one uses a different class of modules as building
blocks.  The techniques for studying general modules in terms of the
simpler ones is via a resolution.

In this paper we propose the class of ``quasi-free Hilbert modules'' as
the building blocks for the general ones.
We discuss the Sz.-Nagy--Foias
model \cite{Na-Fo} in terms of resolutions recalling that this 
interpretation was  a
key motivation for the module approach to multivariate operator
theory.  Then we show under reasonably general hypotheses, 
involving the support of the module in some sense, the existence of a
topologically exact resolution by quasi-free Hilbert modules.  We
discuss various characterizations of the class of modules possessing
such resolutions.  Finally, we describe some
uses of resolutions, that is, how one can extract information and
invariants for a Hilbert module from resolutions.

\subsection{\sc Hilbert Modules}
\noindent
Let $\Omega$ be a bounded domain in $\C^m$.  Examples are the unit
ball, $\B^m$, the polydisk, $\D^m$, or any strongly pseudo-convex
domain $\Omega$ in $\C^m$.  Of course, there are also many examples
for which the boundary, $\partial\Omega$, is not nice.  Nonetheless,
we will consider the natural function algebra $\mathcal{A}(\Omega)$ obtained
from the closure in the supremum norm on $\Omega$ of all functions
holomoprhic in some neighborhood of the closure of $\Omega$.  If
$\partial\Omega$ is not nice, there may be other natural algebras, 
perhaps generated by the polynomials or rational functions with poles
outside the closure of $\Omega$. 
For more refined results, one will probably need to make
additional assumptions about the boundary but we will not need to do
that in this paper.  For $\Omega =$ $\B^m$ or $\D^m$, one obtains the
familiar ball algebra or the poly-disk algebra.

A Hilbert module $\mathcal{M}$ over $\mathcal{A}(\Omega)$ is a Hilbert
space with a multiplication $\mathcal{A}(\Omega) \times \mathcal{M} \to
\mathcal{M}$ making $\mathcal{M}$ into a unital module over
$\mathcal{A}(\Omega)$ and such that multiplication is continuous.
Using the closed graph theorem one can show the existence of a
constant $\alpha$ such that
$$
\|\varphi f \|_\mathcal{M} \leq \alpha \|\varphi\|_{\mathcal{A}(\Omega)}
\|f\|_\mathcal{M}.
$$

One says that $\mathcal{M}$ is a contractive Hilbert module if $\alpha =
1$.  Classical examples of contractive Hilbert modules are 
\begin{enumerate}
\item[\rm (i)] the Hardy module $H^2(\D^m)$ (over the algebra 
$\mathcal{A}(\D^m)$) which is the closure of the polynomials, 
$\C[\ul{z}]$, in $L^2(\partial D^m)$ and    
\item[\rm (ii)] the Bergman module, $B^2(\Omega)$ (over the algebra 
$\mathcal{A}(\Omega)$) which is the closure of $\mathcal{A}(\Omega)$ 
in $L^2(\Omega)$ with volume measure on $\Omega$.
\end{enumerate}

The simplest family of modules over $\mathcal{A}(\Omega)$ corresponds
to evaluation at a point in the closure of $\Omega$.  For $\ul{z}$ in
the closure of $\Omega$, we make the one-dimensional Hilbert space
$\C$ into the Hilbert module $\C_{\ul{z}}$, by setting $\varphi v =
\varphi(z) v$ for $\varphi$ $\in \mathcal{A}(\Omega)$ and $v \in \C$.
If $\mathcal{M}$ and $\mathcal{N}$ are Hilbert modules over
$\mathcal{A}(\Omega)$, then there are two obvious ways to make the
Hilbert space tensor product $\mathcal{M} \otimes \mathcal{N}$ into a
module over $\mathcal{A}(\Omega)$.  One obtains the module tensor
product, $\mathcal{M}\otimes_{\mathcal{A}(\Omega)} \mathcal{N}$, by
requiring that the multiplication on $\mathcal{M}$ and
the multiplication on $\mathcal{N}$ are equal, that is, one takes the
closed submodule $\mathcal{S}$ of $\mathcal{M}\otimes \mathcal{N}$
generated by the elements $\varphi f \otimes g - f \otimes \varphi g$,
for $f \in \mathcal{M}$, $g \in \mathcal{N}$, and $\varphi\in
\mathcal{A}(\Omega)$ and defines
$\mathcal{M}\otimes_{\mathcal{A}(\Omega)} \mathcal{N}$ to be the
quotient $\mathcal{M}\otimes\mathcal{N}/\mathcal{S}$ on which the two
multiplications agree.  

We use the module tensor product to accomplish localization of a
Hilbert module $\mathcal{M}$ over $\mathcal{A}(\Omega)$ by considering
$\mathcal{M} \otimes_{\mathcal{A}(\Omega)} \C_{\ul{z}}$ which is
isomorphic to the direct sum of $k$ copies of $\C_{\ul{z}}$, where $k$
can be any cardinal number from zero to the module rank of
$\mathcal{M}$.  Not only do we consider these tensor products at each
point of the closure of $\Omega$, but together they form the spectral
sheaf, ${\rm Sp}(\mathcal{M})$, of the module $\mathcal{M}$.
We will let ${\rm Sp}_{\ul{z}}(\mathcal{M})$ denote 
$\mathcal{M} \otimes_{\mathcal{A}(\Omega)} \C_{\ul{z}}$. 
The spectral sheaf ${\rm Sp}(\mathcal{M})$ can be the zero sheaf
but, in general, it consists of the direct sum of a number 
of copies of $\C_{\ul{z}}$
over each point $\ul{z}$ in the closure of $\Omega$.  If the spectral sheaf 
${\rm Sp}(\mathcal{M})$ defines a holomorphic hermitian bundle over 
$\Omega$, then we say that $\mathcal{M}$ is {\em locally free}.
The spectral sheaf of the Hardy or Bergman modules is a hermitian
holomorphic line bundle over $\Omega$.

There have been two main lines of research concerning Hilbert modules,
one studying submodules and the other quotient modules.  Using the
classical theorem of Beurling \cite{AB} on invariant subspaces of the
unilateral shift operator, one can show that all nontrivial submodules
of the Hardy module $H^2(\D)$ over the disk algebra $\mathcal{A}(\D)$
are isometrically isomorphic.  The Rigidity Theorem 
(\cite{rgdvip}, \cite{D-P-S-Y})
shows that the situation in higher dimensions is very different.  Two
submodules defined by taking the closure of ideals in $\C[\ul{z}]$ in
the Hardy or Bergman modules (and other more general modules), which
satisfy certain hypotheses, are similar or even quasi-similar if and
only if the ideals are the same.  Thus the rigidity, the closures cannot 
be equivalent in any reasonable sense unless the ideals
are equal.  The hypotheses eliminate principal ideals and insure that
the zero sets of the associated primary ideals intersect $\Omega$.
The result demonstrates that there is a great variety of non
equivalent Hilbert modules in the higher dimensional setting.  The
proof relies on a higher order generalization of the spectral sheaf
and rests on results from commutative algebra.

The work on quotient modules concerns relating properties of a quotient 
module with those of the submodule in cases where the latter consists of 
functions that vanish to some order in the normal direction to a 
hypersurface.  Again subject to mild hypotheses, one characterizes 
(cf. \cite{DM}, \cite{DMC}) the quotient module in terms of the local 
geometry of the hypersurface and the spectral sheaf of the larger module.

\subsection{\sc Sz.-Nagy--Foias Model}
\noindent
One powerful approach to the study of contraction operators on a
complex Hilbert space is the model theory of Sz.-Nagy and Foias
\cite{Na-Fo}.  To understand the interpretation of their model in the
module context, we must first recall the
theorem of von Neumann \cite{von1} which states for a contraction operator
$T$ on a Hilbert space $\mathcal{H}$ and a polynomial $p$ we have
$\|p(T)\| \leq \|p \|_{\mathcal{A}(D)}$.  This inequality enables one
to make $\mathcal{H}$ into a contractive Hilbert module over
${\mathcal{A}(D)}$.  Thus, contraction operators on Hilbert space and
contractive Hilbert modules are two ways of looking at the same thing.

The co-isometric form of the Sz.-Nagy--Foias model for contraction
operators yields an isometry $W$ on a Hilbert space
$\mathcal{K}=\mathcal{H} \oplus \mathcal{G}$ with $W(\mathcal{G})
\subseteq \mathcal{G}$ and such that $T = P_\mathcal{H}
W_{|\mathcal{H}}$.  This can also be written $0 \leftarrow \mathcal{H}
\leftarrow \mathcal{K} \leftarrow \mathcal{G} \leftarrow 0$, where the
arrows are module maps with the map from $\mathcal{K}$ to
$\mathcal{H}$ being the orthogonal projection and the map from
$\mathcal{G}$ to $\mathcal{K}$ being inclusion.  This is an example of
a resolution of Hilbert modules.  Moreover, one has that $\mathcal{G}$
is unitarily equivalent to some vector-valued Hardy module
$H^2_{\mathcal{E}^*}(\D)$ and $\mathcal{K}$ is unitarily equivalent to
$H^2_{\mathcal{E}}(\D) \oplus \mathcal{U}$, where $W$ is the shift
operator on $H^2_{\mathcal{E}}(\D)$ and is a unitary operator on
$\mathcal{U}$.  Further, one shows that $\mathcal{U} = \{0\}$ if and
only if $T^{*k} \to 0$ in the strong operator topology.  Such
contractions are said to belong to class $C_{.0}$ by Sz.-Nagy and
Foias and in this case the resolution has the simpler form $0
\leftarrow \mathcal{H} \leftarrow H^2_{\mathcal{E}}(\D) \leftarrow
H^2_{\mathcal{E}^*}(\D) \leftarrow 0$.    
The modules appearing in such a resolution are the direct sum of
copies of the Hardy module and all the module maps are partial
isometries which makes the resolution especially nice.

No analogous results are known for the bounded case unless the module
is completely bounded and hence similar to a contractive one \cite{vip}.
Here we are interested in the question of when resolutions exist, not
just over $\mathcal{A}(\D)$ but for the multivariate case which we 
take up in the next section.

\section{Quasi-Free Modules}

As we indicated in the introduction, a most important issue in
considering resolutions of Hilbert modules is just what collection of
modules to use as the building blocks.  A second
issue concerns the nature of the module maps.  In the case of class
$C_{.0}$ contractive Hilbert modules over the disk algebra, the modules
used are direct sums of copies of the Hardy module and the maps are
partial isometries.  Moreover, the existence of such a resolution is
based on the existence of a unitary or co-isometric dilation.  Most of
the early consideration of resolutions (\cite{RGD}, \cite{RGD2}) followed this lead
and, for example, the notion of Silov module was introduced for this
reason.  Now, however, constructing resolutions via such dilations
seems to be the wrong approach\footnote{In \cite{WBA4}, 
Arveson reaches the same conclusion, but he shows that dilations
of a different nature seem to work well in his context.}, at least
for Hilbert modules over $\mathcal{A}(\Omega)$, when $\Omega$ lies in 
$\C^m$ with $m > 1$.  To illustrate, we consider a simple example.

If we consider $\C_{\ul{0}}$ over $\mathcal{A}(\D^2)$, then we have that 
$\C_{\ul{0}}$ 
is unitarily equivalent to $H^2(\D^2)/H^2_{\ul{0}}(\D^2)$, where 
$H^2_{\ul{0}}(\D^2) = \{f \in H^2(\D^2) : f(\ul{0})= 0\}$.  
Hence, $H^2(\D^2)$ provides a co-isometric 
dilation of $\C_{\ul{0}}$ and $L^2(\partial \D^2)$ is a unitary dilation of 
$\C_{\ul{0}}$ a la Ando \cite{TA}.  However, $H^2(\D^2)$
and $H^2_{\ul{0}}(\D^2)$ are not unitarily equivalent.  More important, if we
consider their spectral sheaves, then

$${\rm Sp}_{\ul{z}} (H^2(\D^2))
= H^2(\mathbb{D}^2)\otimes_{\mathcal{A}}{\mathbb{C}}_{\underline{z}}\cong 
\mathbb{C}_{\underline{z}},~\underline{z}\in \mathbb{D}^2$$
and 
$${\rm Sp}_{\ul{z}} (H^2_{\ul{0}}(\D^2)) = 
H_{\underline{0}}^2(\mathbb{D}^2)\otimes_{\mathcal{A}}\mathbb{C}_{\underline{z}} \cong
\left\{ \begin{array}{ll}
\mathbb{C}_{\underline{z}} & \underline{z}\neq\underline{0} \\ 
\mathbb{C}_{\underline{0}}\oplus\mathbb{C}_{\underline{0}} & \underline{z}=
\underline{0} \end{array}\right.$$\\

This shows, in particular, that $H^2(\D^2)$ and $H^2_{\ul{0}}(\D^2)$ 
are not similar but also that ${\rm Sp}(H^2_{\ul{0}}(\D^2))$ is not 
a vector bundle.   Therefore, $H^2_{\ul{0}}(\D^2)$ is not locally free.  
However, the resolution:
$$
0 \longleftarrow \mathbb{C}_{\underline{0}} \longleftarrow
H^2(\mathbb{D}^2)\stackrel{X}{\longleftarrow} H^2(\mathbb{D}^2)\oplus 
H^2(\mathbb{D}^2)\stackrel{Y}{\longleftarrow} H^2(\mathbb{D}^2) \longleftarrow 0,
$$
where  $(X (f \oplus g))(\ul{z})= z_2f(\ul{z}) - z_1g(\ul{z})$ 
and  $Y(f)(\ul{z})  =  (z_1 f\oplus z_2 f)(\ul{z})$, demonstrates that 
$\C_{\ul{0}}$ has a locally free resolution. 

Now, we introduce the notion of a quasi-free Hilbert module which will be
the ``nice modules'' we will use for building blocks.  This concept is a
refinement of the notions of sharp kernel and kernel Hilbert space
introduced by Curto and Salinas \cite{rcns}, Agrawal and Salinas \cite{AS}.

Let $k$, $1 \leq k \leq \infty$, be an arbitrary cardinal number
and $\ell^2_k$ denote the $k$ - dimensional Hilbert space. 
Let $\mathcal{M}$ be a Hilbert module of rank $k$ over the algebra
$\mathcal{A}(\Omega)$ relative to the generating set $\{f_1, f_2, 
\ldots \} \subseteq \mathcal{A}(\Omega)$ for which 
$\{f_i \otimes_{\mathcal{A}} 1_{\ul{z}}: 1\leq i \leq k\}$ 
is linearly independent for $\ul{z}\in \Omega$.  Since
module multiplication by a function $\varphi$ on the module tensor
product $\mathcal{M}\otimes_{\mathcal{A}(\Omega)}\C_{\ul{z}}$ is just
multiplication by $\varphi(\ul{z})$, we see it must be isomorphic to
the Hilbert space tensor product $\C_{\ul{z}} \otimes \ell_k^2$.
Define the map $\Gamma: \mathcal{A}(\Omega)\otimes_{\rm alg} \ell^2_k \to
\mathcal{M}$ by $\Gamma( \sum \varphi_i \otimes e_i) = \sum \varphi_i f_i$,
where $\{e_i\}_{i=1}^k$ is the standard orthonormal basis in
$\ell^2_k$.  We claim that $\Gamma$ is well defined, one-to-one and has
dense range.  Given the uniqueness of expressing an element $\phi =
\sum \varphi_i \otimes e_i$ as a finite sum, we have that $\Gamma$
is well defined. If $\Gamma(\phi) = 0$, then we have for $\ul{z} \in
\Omega$ that $\sum \varphi_i(\ul{z})\big ( f_i
\otimes_{\mathcal{A}}1_{\ul{z}}\big ) = \big (\sum \varphi_i f_i \big )
\otimes_{\mathcal{A}}1_{\ul{z}} = \Gamma(\phi) 
\otimes_{\mathcal{A}}1_{\ul{z}} = 0$.  
Since the $f_i$ are linearly independent, it
follows that each $\varphi_i = 0$ and hence $\phi = 0$.  Finally, the
range of $\Gamma$ is dense since the $f_i$'s form a generating set for
$\mathcal{M}$.  Now define the inner product $\langle\ , \rangle$ on
$\mathcal{A}(\Omega) \otimes_{\rm alg} \ell^2_k$ such that $\langle
\phi, \psi \rangle = \langle \Gamma \phi, \Gamma \psi \rangle_\mathcal{M}$.  We
complete $\mathcal{A}(\Omega) \otimes_{\rm alg} \ell^2_k$ using this
inner product to obtain a Hilbert module isometrically isomorphic to
$\mathcal{M}$. Let $e_{\ul{z}}: \mathcal{A}(\Omega) \otimes_{\rm alg} \ell^2_k 
\to  \ell^2_k$ be the evaluation map at $\ul{z} \in \Omega$.

Let $X_{\ul{z}}:\mathcal{M} \otimes_\mathcal{A} \C_{\ul{z}} \to
\C_{\ul{z}}\otimes \ell^2_k$ be the map defined by $X_{\ul{z}}(f_i
\otimes_\mathcal{A} 1_{\ul{z}}) = 1_{\ul{z}} \otimes e_i$, and extend
linearly to finite linear sums.  

\begin{Lemma} \label{eval}
The map $X_{\ul{z}}$ is bounded if the evaluation map $e_{\ul{z}}$ is bounded.  
Furthermore, $\|X_{\ul{z}}\| = \| e_{\ul{z}}\|$, for $\ul{z}\in \Omega$.
\end{Lemma}
\begin{proof} Let $\phi = \sum \varphi_i e_i$ be any arbitrary 
element of $\mathcal{A}(\Omega) \otimes_{\rm alg} \ell^2_k$. 
In the following, $\|\phi\|_{\mathcal{M}}$ denotes the 
norm induced on $\mathcal{A}(\Omega)\otimes_{\rm alg} \ell^2_k$ 
via the map $\Gamma$.
First, let us compute the norm of the operator 
$e_{\ul{z}} : \mathcal{A}(\Omega)\otimes_{\rm alg} \ell^2_k \to \ell^2_k$ 
as follows:
\begin{eqnarray*} \label{Xzez}
\|e_{\ul{z}}\| &=& \sup_\phi \| \phi (\ul{z})\|_{\ell^2_k} / 
\| \phi \|_\mathcal{M} \nonumber\\ 
&=& \sup_\phi \| \sum \varphi_i (\ul{z}) e_i \|_{\ell^2_k} / 
\| \Gamma \phi \|_\mathcal{M} \nonumber\\
&=& \sup_\phi \| \sum \varphi_i(\ul{z}) e_i \|_{\ell^2_k} / 
\|\sum \varphi_i f_i \|_\mathcal{M} \nonumber\\
&=& \sup_{\phi, \Psi} \| \sum \varphi_i(\ul{z}) e_i \|_{\ell^2_k} / 
\|\sum (\varphi_i + \psi_i) f_i \|_\mathcal{M} \nonumber\\
&=& \sup_{\phi} \| \sum \varphi_i(\ul{z}) e_i \|_{\ell^2_k} / 
\inf_{\Psi} \|\sum (\varphi_i + \psi_i) f_i \|_\mathcal{M},
\end{eqnarray*}
where $\Psi=\{\psi_1, \ldots, \psi_n\}$ is a set of finitely many 
non-zero elements in 
$\mathcal{A}(\Omega)$ that vanish at $\ul{z}\in \Omega$.
Note that, $\sum_{1\leq i \leq n} \psi_i f_i$ is a collection of elements dense 
in $\mathcal{A}_{\ul{z}}\mathcal{M}$, where $\mathcal{A}_{\ul{z}}$ is the ideal
of functions in $\mathcal{A}(\Omega)$ that vanish at $\ul{z}$. 
Therefore, we see that 
$$\inf_{\Psi} \|\sum (\varphi_i + \psi_i) f_i \|_\mathcal{M} = 
\| (\sum \varphi_i f_i) \otimes_{\mathcal{A}} 1_{\ul{z}}\|
_{\mathcal{M} \otimes_\mathcal{A} \mathbb{C}{\ul{z}}}.$$
Consequently,  $\|e_{\ul{z}}\| = 
\sup \| \sum_\phi \varphi_i (\ul{z}) e_i \|_{\ell^2_k} / 
\|( \sum \varphi_i f_i)\otimes_{\mathcal{A}} 1_{\ul{z}} \|
_{\mathcal{M} \otimes_\mathcal{A} \mathbb{C}{\ul{z}}}$.
Since $\| \sum \varphi_i (\ul{z}) e_i \|_{\ell^2_k }= 
\| (\sum \varphi_i(\ul{z}) e_i) 
\otimes 1_{\ul{z}}\|_{{\ell^2_k} \otimes \mathbb{C}_{\ul{z}}}$ 
is by definition the norm of 
the operator $X_{\ul{z}}$, it follows that $\|e_{\ul{z}}\| 
= \| X_{\ul{z}}\|$ as claimed.
\end{proof}
This Lemma prompts the following Definition.
\begin{Definition} \label{defqf} 
A Hilbert module $\mathcal{R}$ over $\mathcal{A}(\Omega)$ 
is said to be quasi-free of rank $k$ relative to the generating set 
$\{f_1, f_2,\ldots \}$ for $1 \leq k \leq \infty,$ if
\begin{enumerate}
\item[\rm (i)] $f_1 \otimes_{\mathcal{A}} 1_{\ul{z}}, f_2 
\otimes_{\mathcal{A}} 1_{\ul{z}}, \ldots $ forms a basis for the fiber 
${\rm Sp}_{\ul{z}}(\mathcal{R})$ for $\ul{z} \in \Omega$,
\item[\rm (ii)] the map $X_{\ul{z}}$ is locally uniformly bounded in norm, and 
\item[\rm (iii)]  for $f$ in $\mathcal{R}$, 
$f \otimes_{\mathcal{A}} 1_{\ul{z}} = 0$ for every $\ul{z} \in \Omega$ 
if and only if $f = 0$ in $\mathcal{R}$.
\end{enumerate}
\end{Definition}
When $k$ is finite, the combination of the requirements that 
${\rm Sp}_{\ul{z}}(\mathcal{R})$ is
$k$ - dimensional and the localization of the generating set has
cardinality $k$ has strong implications.  For $k=\infty$, there are many
different ways in which a set can be a basis.  Clearly, we don't want
to assume the set forms an orthonormal or even an orthogonal basis.
But we may want to assume that the set of vectors $\{f_1
\otimes_\mathcal{A} 1_{\ul{z}}, f_2 \otimes_{\mathcal{A}} 1_{\ul{z}},
\cdots \}$ forms a basis in ${\rm Sp}_{\ul{z}}(\mathcal{R})$
equivalent to the standard basis.  However, in this paper we assume
only that the set is linearly independent and spans.

There is another description which 
demonstrates the sense in which a
quasi-free Hilbert module is ``almost free''. In commutative algebra, the
statement that a module is free means that it is isomorphic to a
direct sum of copies of the ring which in our case would be
$\mathcal{A}(\Omega)$.  But this can't happen, since the direct sum of
copies of $\mathcal{A}(\Omega)$ can't be isomorphic to a Hilbert space
unless $\mathcal{A}(\Omega) = \C$ which is impossible.  Hence, we do
the next best thing.  Since we are interested in modules with a
Hilbert space structure, we will begin with the ``free module'' 
$\mathcal{A}(\Omega) \otimes_{\rm alg}\ell^2_k$ which is the
algebraic tensor product of $\mathcal{A}(\Omega)$ with the Hilbert space 
$\ell^2_k$, $1\leq k \leq \infty$, and then complete it to obtain a  
Hilbert space.

A  module $\mathcal{R}$ over $\mathcal{A}(\Omega)$, 
quasi--free relative to $\{ f_1, f_2, \ldots \}$, is the Hilbert space
completion of the free module $\mathcal{A}(\Omega) \otimes_{\rm alg}
\ell^2_k$ via the map $\Gamma$.  Moreover,  the following statement 
is an abstract characterization of such completions of 
$\mathcal{A}(\Omega) \otimes_{\rm alg} \ell^2_k$. 
The technique used in its proof is closely related to the proof of 
\cite[Theorem 5.14]{rgdvip}.

\begin{Proposition} \label{charqf} 
A Hilbert module $\mathcal{R}$ for $\mathcal{A}(\Omega)$
is quasi-free of rank $k$ for $1\leq k \leq \infty$, relative to 
a generating set $\{f_1, f_2, \ldots  \}$ if and only if it is 
isometrically isomorphic to the completion of 
$\mathcal{A}(\Omega) \otimes_{\rm alg} \ell^2_k$ with
respect to a norm associated with an inner product such that 
\begin{enumerate}
\item[\rm (a)] evaluation of functions in
$\mathcal{A}(\Omega) \otimes_{\rm alg} \ell^2_k$ at each point
$\ul{z}$ in $\Omega$ is locally uniformly bounded, 
\item[\rm (b)] module multiplication on 
$\mathcal{A}(\Omega) \otimes_{\rm alg} \ell^2_k$ by
functions in $\mathcal{A}(\Omega)$ is continuous, and 
\item[\rm (c)] for $\{\phi_n\}$ contained in 
$\mathcal{A}(\Omega) \otimes_{\rm alg}\ell_k^2$ 
which is Cauchy in norm, we have $\|\phi_n(\ul{z})\|_{\ell^2_k} \to 0$ 
for all $\ul{z} \in \Omega$ if and only if $\|\phi_n\| \to 0$.
\end{enumerate}
\end{Proposition}
\begin{proof} We first show that the inner product 
introduced on $\mathcal{A}(\Omega)\otimes \ell^2_k$
using the map $\Gamma$ satisfies conditions (a), (b) and (c).
Let $\mathcal{R}$ be a quasi-free module.
Then the local uniform boundedness of the map $X_{\ul{z}}$ 
together with the equality $\|X_{\ul{z}}\| = \|e_{\ul{z}}\|$, 
for $\ul{z} \in \Omega$ establishes the condition (a).
For the proof of (b), consider the function $\psi$ in 
$\mathcal{A}(\Omega)$ and observe that $\|\psi \sum \varphi_i
\otimes e_i\| = \|\sum \psi \varphi_i \otimes e_i\| = 
\|\Gamma (\sum \psi \varphi_i \otimes e_i)\|_{\mathcal{R}} = 
\|\psi \Gamma(\sum \varphi_i \otimes e_i)\|_{\mathcal{R}} \leq C_1 
\|\psi\| \|\Gamma(\sum \varphi \otimes e_i)\|_{\mathcal{R}} \leq C_1 
\|\psi\| \|\sum (\varphi_i \otimes e_i)\|$.  
Finally, let $\phi_n = \sum \varphi_i^{(n)} \otimes e_i$ be a
sequence in $\mathcal{A}(\Omega) \otimes_{\rm alg} \ell^2_k$, which is
Cauchy in norm.  Then $\Gamma(\phi_n)\to g$ for some $g \in
\mathcal{R}$.  Since $\Gamma$ is continuous, it follows that
$e_{\ul{z}}(\phi_n) \to 0$ if and only if
$e_{\ul{z}}(\Gamma(\phi_n))\to 0$.  Or, in other words,
$\phi_n(\ul{z}) \to 0$ if and only if $\sum \varphi_n(\ul{z}) f_i
\otimes_{\mathcal{A}} 1_{\ul{z}}\to 0$.  Hence, the assumption that
$\phi_n(\ul{z})\to 0$ implies that $g(\ul{z}) = 0$ and hence $g = 0$
by (iii) of Definition \ref{defqf}.  This shows that the condition (c) holds
which completes the proof in the first direction.
  
For the converse, assume that $\mathcal{R}$ is the completion of
$\mathcal{A}(\Omega) \otimes_{\rm alg} \ell^2_k$ with respect to an
inner product that satisfies (a), (b), and (c) of the Proposition.  We
must verify that the conditions of Defintion 1 hold.  Fix $\ul{z}\in
\Omega$ and consider the map $F_{\ul{z}}:\mathcal{A}(\Omega)
\otimes_{\rm alg}\ell^2_k \to \ell^2_k$ defined by $F_{\ul{z}}(\sum \varphi_i
\otimes e_i) = \sum \varphi_i(\ul{z}) e_i$.  By condition (a) of the
Proposition, it follows that $F$ extends to a map from $\mathcal{R}$
to $\ell^2_k$.  We use $F_{\ul{z}}$ to define a map $F^\prime_{\ul{z}}: \mathcal{R}
\otimes \C_{\ul{z}} \to \C_{\ul{z}} \otimes \ell^2_k$ such that
$F^\prime_{\ul{z}}\big ((\sum \varphi_i \otimes e_i \big ) \otimes a = \sum
\varphi_i(\ul{z}) e_i \otimes a$. The kernel of this map 
is spanned by the vectors $\varphi \otimes
e_i \otimes a - 1 \otimes e_i \otimes \varphi(\ul{z})a$ for $\varphi_i \in
\mathcal{A}(\Omega)$ and $a \in \C_{\ul{z}}$, which is the submodule
used to define the module tensor product
$\mathcal{R}\otimes_{\mathcal{A}}\C_{\ul{z}}$.  Hence, we see that
evaluation of $\sum \varphi_i \otimes e_i$ at $\ul{z}$ matches $\sum
\varphi_i (\ul{z}) e_i$ in $\C_{\ul{z}} \otimes \ell^2_k$ which
establishes (i) in Definition \ref{defqf}.  The condition (ii) of the Definition 
is clearly the same as condition (a) of the Proposition.  Also, this matchup shows
that condition (c) of the proposition implies (iii) of Defintiion 1,
which completes the proof.
\end{proof}

Observe that no assumption of holomorphicity is made in the definition
of a quasi-free Hilbert module $\mathcal{R}$.  However, identification
of $\mathcal{R}$ with the completion of $\mathcal{A}(\Omega)
\otimes_{\rm alg} \ell^2_k$ in the finite rank case 
makes the spectral sheaf ${\rm Sp}(\mathcal{R})$
into a hermitian holomorphic vector bundle of rank $k$ with 
holomorphic frame 
$\ul{z} \to \{f_i \otimes_{\mathcal{A}}1_{\ul{z}}: 1\leq i \leq k\}$. 
Moreover, using the identification of $\mathcal{R}$ with the
completion of $\mathcal{A}(\Omega) \otimes_{\rm alg} \ell^2_k$, we see
that $\mathcal{R}$ can be realized as a space of $\ell^2_k$ - valued
holomorphic functions on $\Omega$ which forms a kernel Hilbert space.

Obviously, the Hardy and Bergman modules are quasi-free.  While 
submodules of quasi-free modules are not quasi-free in general, 
principal submodules are, since one can view them as being obtained 
merely from a change of norm.  Quotient modules of quasi-free Hilbert 
modules are seldom quasi-free.  However, it can happen.  The 
relationship here is analogous to the situation of holomorphic 
subbundles of holomorphic bundles; sometimes there is a holomorphic 
complement.   The following statement should be true in our context and 
would clarify the relation between the notions of free and projective.  

\begin{Conjecture}
  Let $\mathcal{R}$ be a quasi-free Hilbert module of rank $k$, $1
  \leq k < \infty$, over $\mathcal{A}(\Omega)$ and $\mathcal{R}_1$ and
  $\mathcal{R}_2$ be submodules of $\mathcal{R}$ such that
  $\mathcal{R}$ is the algebraic direct sum of $\mathcal{R}_1$ and
  $\mathcal{R}_2$.  Then $\mathcal{R}_1$ and $\mathcal{R}_2$ are
  quasi-free of ranks $k_1$ and $k_2$, respectively, and $k = k_1 +
  k_2$.
\end{Conjecture}

Something analogous should be true in the case $k = \infty$ but would
probably require more explicit hypotheses on the angle between the two
submodules. 

\section{Regular Modules}

As indicated earlier, a resolution of the Hilbert module $\mathcal{M}$ 
is a sequence of modules $\{\mathcal{R}_i\}$, either of finite or 
infinite length, with module map $X_0: \mathcal{R}_1 \to \mathcal{M}$ 
and module maps $X_i : \mathcal{R}_{i+1} \to \mathcal{R}_i \mbox{~for~}i \geq 1$, 
such that the sequence  
$$                                                                
0 \longleftarrow \mathcal{M} \stackrel{X_0}{\longleftarrow}\mathcal{R}_1
\stackrel{X_1}{\longleftarrow} \mathcal{R}_2 \stackrel{X_2}{\longleftarrow}\cdots
$$
is exact in the sense that $X_0$ is onto and  $\ker X_i = \mbox{ran}\, X_{i+1}$ for 
$i\geq 1$.  If the sequence is of finite length, then we must have the final 
$\mathcal{R}_i = 0$.  One speaks of a weak resolution if it is only 
topologically exact or, equivalently, if one assumes that 
$X_0$ has dense range and $\ker X_i = \mbox{clos ran}\, X_{i+1}$ for $i \geq 1$.  
We are assuming in all cases that the modules 
$\{\mathcal{R}_i\}$ are quasi-free over $\mathcal{A}(\Omega)$.  
One can also put an additional 
restriction on the module maps by requiring that they are partial 
isometries in which case we will speak of a strong resolution.

The resolution obtained from the Sz.-Nagy--Foias model is a strong 
resolution, while the second one constructed for the Hilbert module 
$\C_{\ul{0}}$ over $\mathcal{A}(\D^2)$ is a resolution but not a 
strong one. Although one seeks 
resolutions as nice as possible, and closely related to $\mathcal{M}$, 
one often faces tradeoffs.  For example, the inclusion map of the Hardy module 
$H^2(\D)$ into the Bergman module $B^2(\D)$ defines a weak resolution with just 
one term since the map has dense range and trivial kernel.  However, it 
is not clear just what this resolution can tell us about the Bergman 
module in terms of the Hardy module.  On the other hand, while the 
resolution of the Bergman module given by the Sz.-Nagy--Foias model is a 
strong one, it is obtained at the price of having to deal with  
Hardy modules of infinite multiplicity.  Still we show in section 5 that 
all resolutions, even weak ones, contain information about the 
module.  Finally, in that section we will also compare the existence 
questions for the various kinds of resolutions.  But first we want to 
investigate existence in general.

If one is given a Hilbert module $\mathcal{M}$, the first task in
constructing a resolution of $\mathcal{M}$ by quasi-free Hilbert
modules is to obtain a quasi-free Hilbert module $\mathcal{R}$ and
module map $X: \mathcal{R} \rightarrow \mathcal{M}$ with dense range.
We introduce the following definitions to capture the kinds of
behavior possible with regard to the construction of resolutions.
\begin{Definition} \label{defreg}
A Hilbert module  over $\mathcal{A}(\Omega)$ is said to be
\begin{enumerate}
    \item weakly regular if there exists a quasi-free Hilbert module $\mathcal{R}$ 
over $\mathcal{A}(\Omega)$ and a module map  $X: \mathcal{R} \rightarrow \mathcal{M}$ 
with dense range,
    \item  regular if the map $X$ can be taken to be onto,
    \item  strongly regular if the map $X$ can be taken to be a co-isometry, and
    \item  singular if the only module map $X: \mathcal{R} \rightarrow \mathcal{M}$ 
from a quasi-free Hilbert module $\mathcal{R}$ to $\mathcal{M}$ is the zero map.
\end{enumerate}
\end{Definition}

It is relatively straight forward to see that not 
all Hilbert modules are weakly regular.
In particular, we will see that for the Hilbert module $\C_1$ over 
$\mathcal{A}(\D)$, the only module map $X: \mathcal{R} \to \C_1$ for 
$\mathcal{R}$ a quasi-free Hilbert module, is the zero map.  
Toward that end, we recall an extension of a notion of
Sz.-Nagy and Foias \cite{Na-Fo} to the context of Hilbert modules 
(cf. \cite{AD}).

\begin{Definition} \label{defC.0}
A Hilbert module $\mathcal{M}$ over $\mathcal{A}(\Omega)$ 
is said to belong to class 
$C_{\cdot 0}$ if for every sequence 
$\{\varphi_n\}_{n\in \mathbb{N}}$ in the unit ball of 
$\mathcal{A}(\Omega)$ satisfying 
$\varphi_n(\ul{z}) \rightarrow 0$ for $\ul{z}$ in $\Omega$ it follow that  
$M_{\varphi_n}^* \rightarrow 0$ in the strong operator topology.
\end{Definition}
One could also assume that the $\varphi_n$ converge uniformly 
to zero on compact subsets of $\Omega$.  In many situations, these two notions 
seem to coincide but it is not clear if they do for general Hilbert 
modules.  
\begin{Proposition} \label{regC.0}
A regular Hilbert module over $\mathcal{A}(\Omega)$ belongs to class 
$C_{\cdot 0}$.
\end{Proposition}

\begin{proof} Let $\mathcal{M}$ be a regular Hilbert module with
$\mathcal{R}$ a quasi-free Hilbert module and $X: \mathcal{R} \to
\mathcal{M}$ an onto module map.  If $k_z$ in $\mathcal{R}$ is a
common eigenvector for the adjoint of module multiplication on
$\mathcal{R}$, then $M_{\varphi_n}^* k_z = \overline{\varphi_n(z)}
k_z$ and hence $\|M_{\varphi_n}^* k_z\| \rightarrow 0$.  Since the
vectors $\{k_z\}$, $z\in \D$ span $\mathcal{R}$ and
$\|M_{\varphi_n}^*\| \leq \alpha \|\varphi_n\| \leq \alpha$, it
follows that $M_{\varphi_n}^* \rightarrow 0$ in the strong operator
topology.  Then, $X^*N_\varphi^* = M_\varphi^*X^*$, where $N_\varphi$
denotes the operator defined by module multiplication on $\mathcal{M}$.
Since $X$ is onto, it follows that $X^*$ is bounded below and hence $N_\varphi^*
\to 0$ in the strong operator topology.  Thus $\mathcal{M}$ belongs to class 
$C_{\cdot 0}$.
\end{proof}

Taking $X$ to be the identity map, this result yields a property of
quasi-free Hilbert modules.

\begin{Corollary} \label{qfC.0}
All quasi-free Hilbert modules belong to class $C_{\cdot 0}$.
\end{Corollary}

\begin{Proposition} \label{C_1}
The Hilbert module $\C_1$ over $\mathcal{A}(\D)$ does not belong to class 
$\C_{\cdot 0}$.
\end{Proposition}

\begin{proof} There exists a sequence $\{\varphi_n\}_{n\in \mathbb{N}} \in
\mathcal{A}(\D)$ satisfying $\varphi_n(1) = 1$, $\|\varphi_n\| = 1$,
and $\varphi_n(z) \rightarrow 0$ for $z \in \D$.  For example, one
can take $\varphi_n(z) = 1/n (1 + 1/n - z)$. Then for $\lambda$ in $\C_1$ we
see that $M_{\varphi_n}^* \lambda = \overline{\varphi_n(1)}\lambda = \lambda$
does not converge to zero which completes the proof.
\end{proof}

We have been unable to show either that weakly regular implies class $C_{\cdot 0}$ 
or that class $C_{\cdot 0}$ implies weak regularity. However, the first 
conclusion would seem to be likely.

Although we have been unable to obtain an intrinsic characterization
of weak regularity, we can provide two properties each of which
implies it, both relating to the notion that the module is supported
on the interior of $\Omega$.

\begin{Definition} \label{smooth}
Let $\mathcal{R}$ be a quasi-free Hilbert module of rank one over
$\mathcal{A}(\Omega)$ for the generating vector $g$ and let $\mathcal{M}$ 
be a Hilbert module over $\mathcal{A}(\Omega)$.
Then $\mathcal{M}$ is said to be smooth relative to $\mathcal{R}$ 
and $g$ if the map 
$S: \mathcal{M} \to \mathcal{M}\otimes_{\mathcal{A}(\Omega)} \mathcal{R}$ 
defined by $Sf = f \otimes_{\mathcal{A}(\Omega)} g$ is one-to-one.
\end{Definition}

Smooth modules in this sense are not always in class $C_{\cdot 0}$.
Consider the contractive Hilbert module over $\mathcal{A}(\D)$ defined
by a function $\theta$ in $H^\infty(\D)$ with $|\theta (e^{it})| < 1$ on
a set of positive measure.  It has the form $\mathcal{H} =
H^2(\D)\oplus L^2(\overline{\Delta})/ \{\theta f \oplus \Delta f : f \in
H^2(\D)\}$, where $\Delta (e^{it}) = (1 - |\theta (e^{it})|^2)^{1/2}$
and $\overline{\Delta}$ is the characteristic function for the support of
$\Delta$.  The map $S: H^2(\D) \oplus L^2(\overline{\Delta} \T) \to
(H^2(\D) \oplus L^2(\overline{\Delta} \T)) \otimes_{\mathcal{A}(\mathbb{D})}
H^2(\D)$ reduces to $S(f \oplus g) = f \otimes_{\mathcal{A}(\mathbb{D})} 1$
since $L^2(\D) \otimes_{\mathcal{A}(\mathbb{D})} H^2(\D) = \{0\}$.  Since the
map from $H^2(\D)$ to $H^2(\D) \otimes_{\mathcal{A}(\mathbb{D})} H^2(\D)$ is
one-to-one, we see that $S(h \oplus k)= 0$ implies $h = 0$.  For such
a vector to be in $\mathcal{H}$, it must be orthogonal to the subspace
$\{\theta f \oplus \Delta f : f \in H^2(\D)\}$.  We can choose $\theta$
so that this subspace is dense in $L^2(\overline{\Delta} \T)$ in which
case $\mathcal{H}$ is $H^2(\D)$-smooth.  One can show that such a
Hilbert module does not belong to the class $C_{\cdot 0}$.

In general, the question of whether $\mathcal{H}$ is $H^2(\D)$-smooth
depends on the density of $\{\Delta f\}$ in $L^2(\overline{\Delta} \T)$ and
that happens when the associated contraction on $\mathcal{H}$ has no
co-isometric part.  That relationship seems to hold for general
contractive Hilbert modules although a precise result would depend
on an analysis of how the notion of smoothness depends on the
quasi-free Hilbert module $\mathcal{R}$ and the generating vector used to define the
$S$-map.  In particular, although we presume the class to be
independent of both $\mathcal{R}$ and $g$, we have been unable to prove that.

\begin{Definition} \label{defSM}
A Hilbert module is said to belong to the class {\em (SM)} if
it is smooth for some quasi-free Hilbert module $\mathcal{R}$ and 
generating vector $g$.
\end{Definition}
Although the notion of smoothness is conceptually appealing, it is not
always easy to verify.  We introduce a subclass, whose membership is
more closely related to operator theoretic ideas which we recall
first.
\begin{Definition} \label{defgeneigen}
If $\mathcal{M}$ is a Hilbert module over $\mathcal{A}(\Omega)$, 
then a vector $h \in \mathcal{M}$ is said to be a common generalized 
eigenvector for the adjoint of module multiplication if $h$ lies in  
the kernel of $(M_\varphi - \varphi(\ul{z})I)^{*n}$  for all $\varphi$ in 
$\mathcal{A}(\Omega)$ and some fixed positive integer $n$.
\end{Definition}
The subclass of the class (SM) we want to define 
involves the assumption that common generalized 
eigenvectors span.
\begin{Definition} \label{defPS}
A Hilbert module is said to belong to the class {\em (PS)} if $\mathcal{M}$ is 
spanned by the generalized eigenvectors for the adjoint of module 
multiplication.
\end{Definition}
Another characterization of class (PS) is possible using the notion of 
higher order localization.  We will not discuss this notion in
complete detail here but see \cite{DMC}.  Consider a point $\ul{z}$ in 
$\Omega$ and let $\mathcal{I}_{\ul{z}}$ be the ideal of polynomials in 
$\C[\ul{z}]$ that vanish at $\ul{z}$ and
$\mathcal{I}_{\ul{z}}^n$ the ideal generated by the products of $n$ 
elements in $\mathcal{I}_{\ul{z}}$.  The
quotient $\mathcal{Q}_{\ul{z}}^n = \C[\ul{z}]/\mathcal{I}_{\ul{z}}^n$ 
with the Hilbert space structure in which
the set of mononials in the quotient form an orthonormal basis, is a
Hilbert module over $\mathcal{A}(\Omega)$ of dimension $mn$ in which module
multiplication by a function $\varphi$ depends on the values at $\ul{z}$ of 
$\varphi$ and its partial derivatives of order less than $n$.  Let $e$ denote the
element of $\mathcal{Q}_{\ul{z}}^n$ corresponding to the monomial $1$.  
Using elementary calculations, one can show that the following proposition provides
another characterization of the class (PS).

\begin{Proposition} \label{charPS} A Hilbert module $\mathcal{M}$ over
$\mathcal{A}(\Omega)$ belongs to class {\em (PS)}
if and only if for every nonzero vector $f \in \mathcal{M}$, there
exists $\ul{z} \in \Omega$ and $n$ such that the image 
$f\otimes_\mathcal{A} e$ in 
$\mathcal{M}\otimes_{\mathcal{A}}\mathcal{Q}_{\ul{z}}^n$ of $f$ is not
$0$.  Equivalently, the intersection of the closures of the submodules
$M(\ul{z}_1,n_1;\ul{z}_2,n_2;\ldots;\ul{z}_k,n_k)$ is $\{0\}$, where 
$M(\ul{z}_1,n_1;\ldots ;\ul{z}_k,n_k)$
is the closure of the product 
$\mathcal{I}_{\ul{z}_1}^{n_1}\cdots \mathcal{I}_{\ul{z}_k}^{n_k}$
acting on $\mathcal{M}$ for
every finite subset ${\ul{z}_1,\ldots ,\ul{z}_k}$ of $\Omega$ and positive integers
${n_1, \ldots ,n_k}$.
\end{Proposition}

\begin{proof}
The equivalence of the two statements is an easy exercise
involving the relation of the submodule which defines the tensor
product of $\mathcal{M}$ with $\mathcal{Q}_{\ul{z}}^n$ and the closure of the range of 
$\mathcal{I}_{\ul{z}}$ acting on $\mathcal{M}$ (cf. \cite[ Theorem 5.14]{rgdvip}).
Similarly, by considering the relation of the kernel of the
adjoint action of the nth power of a function $\varphi$ which vanishes at $\ul{z}$
and the latter space, one shows the equivalence with $\mathcal{M}$ belonging to
class (PS).
\end{proof}

It is not true that the modules belonging to class (PS)  are just those 
determined by their spectral sheaves which are, of course, in class (PS); 
the higher multiplicity examples 
in \cite{DMC} show otherwise.  However, quasi-free Hilbert modules are 
determined by their spectral sheaves and there are relationships 
between the classes $C_{\cdot 0}$ and (PS).

\begin{Corollary} \label{qfPS}
A quasi-free Hilbert module for $\mathcal{A}(\Omega)$ belongs to class {\em (PS)}. 
\end{Corollary}
\begin{Proposition} \label{PSC.0}
A Hilbert module in class {\em (PS)} belongs to class $C_{\cdot 0}$.
\end{Proposition}
The proof is essentially the same as that of Proposition~\ref{regC.0}.  The 
converse to the proposition is false, that is, Hilbert modules of class 
$C_{\cdot 0}$ do not necessarily belong to class (PS).  
For example, let $\theta$ be a singular inner function on the unit disk 
and $\mathcal{M}$ be the quotient Hilbert module  $H^2(\D)/\theta H^2(\D)$.  
Then the quotient map from the quasi-free Hilbert module $H^2(\D)$ onto 
$\mathcal{M}$ shows that it is $C_{\cdot 0}$ but there 
are no common eigenvectors in $\D$ for the adjoint of module multiplication 
since the spectrum of $\mathcal{M}$ is the closed support of the singular measure 
that defines $\theta$ and hence a subset of $\partial \D$.
\begin{Proposition} \label{(SM)}
A Hilbert module in class {\em (PS)} is in class {\em (SM)}.
\end{Proposition}

\begin{proof} Let $\mathcal{M}$ be in class (PS) and let $\mathcal{R}$ 
be a quasi-free Hilbert module of rank one over $\mathcal{A}(\Omega)$ with
generating vector $g$ and $S$ be the map 
$S: \mathcal{M} \to \mathcal{M} \otimes_{\mathcal{A}(\Omega)}
\mathcal{R}$
defined by $Sk = k\otimes_{\mathcal{A}(\Omega)} g$ for $k  \in
\mathcal{M}$.  For $f$ a vector in $\mathcal{M}$, there exists a 
point $\ul{z}\in \Omega$ and an integer $n$ such that 
$f \otimes_{\mathcal{A}(\Omega)} e \not= 0$ in 
$\mathcal{M} \otimes_{\mathcal{A}(\Omega)} \mathcal{Q}$,
where the module $\mathcal{Q}$ stands for $\mathcal{Q}_z^n$ and the 
$e$ as defined above. 
Let $X$ be the map from $\mathcal{M}$ to 
$\mathcal{M}\otimes_{\mathcal{A}(\Omega)}\mathcal{Q}$ defined by 
$X h = h \otimes_{\mathcal{A}(\Omega)}e$ and consider the diagram 
$$
\begin{CD} 
\mathcal{M} @>S_1>>  \mathcal{M} \otimes_{\mathcal{A}(\Omega)} \mathcal{R}\\
@VXVV @VVX\otimes_{\mathcal{A}(\Omega)}I_{\mathcal{Q}}V\\
\mathcal{M} \otimes_{\mathcal{A}(\Omega)} \mathcal{R}@>>S_2> 
\mathcal{M} \otimes_{\mathcal{A}(\Omega)} \mathcal{Q} 
\otimes_{\mathcal{A}(\Omega)} \mathcal{R}\\
\end{CD}
$$
Then, $(X \otimes_{\mathcal{A}(\Omega)}I_{\mathcal{Q}}) S_1 f 
= S_2 X f = S_2 (f \otimes_{\mathcal{A}(\Omega)} e)\not = 0$ 
since $\mathcal{Q} \to \mathcal{Q}\otimes_{\mathcal{A}(\Omega)}
\mathcal{R}$ 
is an isomorphism.  
Hence, $S_1 f \not= 0$ and $\mathcal{M}$ is smooth for $\mathcal{R}$ and $g$ which
completes the proof.
\end{proof}
\begin{Proposition} \label{submodPSC.0}
If $\mathcal{M}$ is a Hilbert module over $\mathcal{A}(\Omega)$ and 
$\mathcal{M}_0$ is a submodule of $\mathcal{M}$, then $\mathcal{M}$ 
in class $C_{\cdot 0}$ implies $\mathcal{M}_0$ is in class $C_{\cdot 0}$ 
and $\mathcal{M}$ in class {\em (PS)} implies $\mathcal{M}_0$ is in class 
{\em (PS)} and $\mathcal{M}$ in class {\em (SM)} implies $\mathcal{M}_0$ 
is in class {\em (SM)}.
\end{Proposition}
These results follow using similar arguments as before.
\begin{Corollary} \label{subqfC.0}
A submodule of a quasi-free Hilbert module over $\mathcal{A}(\Omega)$ belongs 
to classes $C_{\cdot 0}$, {\em (PS)} and {\em (SM)}.
\end{Corollary}
We show in the next section that Hilbert modules of class (PS) and (SM) are 
weakly regular.  With  this result and the above corollary which 
we can use to show that every weakly regular Hilbert module has a 
weak resolution, we will conclude that Hilbert modules in class 
(PS) and (SM) have weak resolutions.  

There is another class of Hilbert modules, which includes those in
class (PS), for which we can also establish the existence of weak
resolutions.

\begin{Definition} \label{defsupp}
A Hilbert module $\mathcal{M}$ over $\mathcal{A}(\Omega)$ is
said to be compactly supported on a vector f in $\mathcal{M}$ if there
is a compact subset $X$ of $\Omega$ and a constant $\beta$ such that 
$$\|\varphi f\|_\mathcal{M} \leq \beta \|\varphi\|_X \|f\|_{\mathcal{M}} 
\mbox{~for~}\varphi \in \mathcal{A}(\Omega),$$
where $\|\varphi\|_X$ denotes the supremum norm of
$\varphi$ taken on $X$.  The module $\mathcal{M}$ is said to be compactly
supported if there is a compact set $X$ and a constant $\beta$ which works for
all vectors $f$ in $\mathcal{M}$.  Finally, the module $\mathcal{M}$ is
said to be almost compactly supported if $\mathcal{M}$ is the span of
the compactly supported vectors in $\mathcal{M}$, where the compact
set and constant can depend on a vector.
\end{Definition}

\section{Construction of Resolution}

We now introduce our basic construction for establishing weak
regularity.  Let $\mathcal{M}$ be a Hilbert module over
$\mathcal{A}(\Omega)$ with a set of $k$ generators $\{f_i \in
\mathcal{M}: 1 \leq i\leq k\}$,  and let
$\mathcal{R}$ be any quasi-free rank $k$ Hilbert module over
$\mathcal{A}(\Omega)$ relative to $\{g_i \in\mathcal{R} : 1 \leq i \leq k \}$, 
$1\leq k \leq \infty$.   
Let $\mathcal{R}_\mathcal{M}$ be the closed submodule of $\mathcal{M}
\oplus \mathcal{R}$ spanned by the vectors
$$  
\{\varphi_i f_i \oplus \varphi_i g_i: \varphi_i \in 
\mathcal{A}(\Omega),\: 1\leq i \leq k \}  
$$
and let $P$ and $Q$ be the module maps from $\mathcal{R}_\mathcal{M}$  
to $\mathcal{M}$ and $\mathcal{R}$, respectively, defined by 
$P(\varphi_i f_i \oplus \varphi_i g_i) = \varphi_i f_i$
and  $Q(\varphi_i f_i \oplus \varphi_i g_i) = \varphi_i g_i$, 
respectively.  Since a dense set of vectors in 
$\mathcal{R}_\mathcal{M}$ has the form  
$\Phi = \sum \varphi_i f_i \oplus \varphi_i g_i$, 
where at most finitely many of the 
$\varphi_i$ are non-zero, we see that  $\|P \Phi\| = 
\|\sum \varphi_i f_i\| \leq (\|\sum \varphi_i f_i\|^2 + 
\|\sum \varphi_i g_i\|^2)^{1/2} 
= \|\Phi\|$ and similarly, $\|Q \Phi\| \leq 
\|\Phi\|$.  Thus both $P$ and $Q$ are well-defined, contractive and have 
dense range.  If $\mathcal{R}_\mathcal{M}$ is quasi-free, 
then it will follow that 
$\mathcal{M}$ is weakly  regular.  

Consider the operator $e_{\ul{z}}: \mathcal{A}(\Omega) \otimes_{\rm alg}
\ell^2_k \to \ell^2_k$ of evaluation at ${\ul{z}}$ in $\Omega$.  
For the function 
$\phi = \sum \varphi \otimes e_i$ in 
$\mathcal{A}(\Omega)\otimes_{\rm alg}\ell^2_k$ we have $e_{\ul{z}}( \phi) = 
\sum \varphi_i(\ul{z}) e_i$, and hence 
\begin{eqnarray*}
\|e_{\ul{z}}\| &=& \sup \big \{ \|e_{\ul{z}}(\phi)\| / \|\phi\|: \phi \in 
\mathcal{A}(\Omega)\otimes_{\rm alg} \ell^2_k \big \},\\
&=& \sup \big \{ \| \sum \varphi_i ({\ul{z}}) e_i\| /
\|\phi\|: \phi \in \mathcal{A}(\Omega)\otimes_{\rm alg} \ell^2_k \big \}.
\end{eqnarray*}
Now, consider the evaluation, first on $\mathcal{R}$ and then on 
$\mathcal{R}_{\mathcal{M}}$, at $\ul{z}$.  Recall that  
$\mathcal{A}(\Omega)\otimes_{\rm alg} \ell^2_k$ is a dense spanning set 
in both $\mathcal{R}$ and $\mathcal{R}_{\mathcal{M}}$. It is clear that 
$$\|e_{\ul{z}}\|_{\mathcal{R} \to \ell^2_k}= 
\sup \big \{ \| \sum \varphi_i ({\ul{z}}) e_i\| /
\|\phi\|: \phi \in \mathcal{A}(\Omega)\otimes_{\rm alg} \ell^2_k \big \}.
$$
However,  
\begin{eqnarray*}
\|e_{\ul{z}}\|_{\mathcal{R}_{\mathcal{M}} \to \ell^2_k} 
&=& \sup \big \{ \|e_{\ul{z}}(\phi )\| / \| \Gamma(\phi)\|: \phi \in 
\mathcal{A}(\Omega)\otimes_{\rm alg} \ell^2_k \big \},\\
&=& \sup \big \{ \| \sum \varphi_i ({\ul{z}}) e_i\| /
\|\sum \varphi_i(f_i \oplus (1\otimes e_i)) \|: 
\phi \in \mathcal{A}(\Omega)\otimes_{\rm alg} \ell^2_k \big \},
\end{eqnarray*}
where $\Gamma: \mathcal{A}(\Omega)\otimes_{\rm alg} \ell^2_k \to 
\mathcal{R}_\mathcal{M}$ is the map defined by 
$\Gamma(\phi) = \sum \varphi_i(f_i \oplus (1\otimes e_i))$. 
So, it follows that the norm of the evaluation operator  on 
$\mathcal{R}$ dominates that of $\mathcal{R}_{\mathcal{M}}$ 
at $\ul{z}$.

\begin{Lemma} \label{nullQ}
The Hilbert module $\mathcal{R}_\mathcal{M}$ over $\mathcal{A}(\Omega)$ is quasi-free 
of rank $k$ relative to the  
generating set $\{f_i \oplus g_i\}$ if $\ker Q = \{0\}$.
\end{Lemma}

\begin{proof} We will establish the three properties for a module to be
quasi-free given in Definition~\ref{defqf}.  
First, the norm of the evaluation operator 
$\|e_{\ul{z}}\|_{\mathcal{R}_\to \ell^2_k}$
is locally uniformly bounded. Since 
$\|e_{\ul{z}}\|_{\mathcal{R}_{\mathcal{M}} \to \ell^2_k} \leq 
\|e_{\ul{z}}\|_{\mathcal{R}_\to \ell^2_k}$, it follows that
the norm of the evaluation operator 
$\|e_{\ul{z}}\|_{\mathcal{R}_{\mathcal{M}}}$ is 
locally uniformly bounded as well, which establishes (ii). 
 
Now, for $\ul{z}$ in
$\Omega$ and $1\leq i \leq k$, we consider the set 
$\{(f_i \oplus g_i) \otimes_\mathcal{A} 1_{\ul{z}}\}$ in 
$\mathcal{R}_\mathcal{M} \otimes_\mathcal{A} \C_{\ul{z}}$.  
Since the set $\{f_i \oplus g_i\}$ generates
$\mathcal{R}_\mathcal{M}$, it follows that the set ${(f_i \oplus g_i)
\otimes_\mathcal{A} 1_{\ul{z}}}$ generates $\mathcal{R}_\mathcal{M}
\otimes_\mathcal{A} \C_{\ul{z}}$ as a module over
$\mathcal{A}(\Omega)$.  However, from the fact that
$\mathcal{R}_\mathcal{M} \otimes_\mathcal{A} \C_{\ul{z}}$ is
isomorphic to a direct sum of copies of $\C_{\ul{z}}$, it follows that
a generating set over $\mathcal{A}(\Omega)$ is the same as a
generating set over $\C$ or as a Hilbert space.  Further, since the
set of vectors $\{g_i \otimes_\mathcal{A} 1_{\ul{z}} \}$ is linearly
independent in $\mathcal{R}\otimes_{\mathcal{A}} \C_{\ul{z}}$, 
it follows that the set $\{(f_i \oplus g_i)\otimes_\mathcal{A}
1_{\ul{z}}\}$ is linearly independent in 
$\mathcal{R}_{\mathcal{M}}\otimes_{\mathcal{A}} \C_{\ul{z}}$, which is
condition (i).    

Thus, whether
$\mathcal{R}_\mathcal{M}$ is quasi-free comes down to whether or not
condition (iii) holds.  Suppose $h$ is a vector in
$\mathcal{R}_\mathcal{M}$ such that $h \otimes_\mathcal{A} 1_{\ul{z}} =
0$ for every $\ul{z} \in \Omega$.
Then we have
$$
(Qh) \otimes_\mathcal{A} 1_{\ul{z}} = 
(Q \otimes_\mathcal{A} 1_{\ul{z}})(h \otimes_\mathcal{A} 1_{\ul{z}}) = 0
$$ 
and, since $Qh$ is in $\mathcal{R}$ which is quasi-free, we have $Qh = 0$.  
If $\ker Q = \{0\}$,  then $\mathcal{R}_\mathcal{M}$ is quasi-free
which completes the proof.
\end{proof}

A reasonable question that arises is whether $\ker Q = \{0\}$ always holds.
To see that is not the case, consider $\mathcal{R} = H^2(\D)$ and 
$\mathcal{M} = \C_1$ over
$\mathcal{A}(\D)$.  One can either repeat the arguments from the last
section or observe that if $\ker Q = \{0\}$ in this case, it would
follow from the lemma and later results in this section that $\C_1$ is
weakly regular.  Since $\C_1$ is finite dimensional, we have that
$\C_1$ is regular which would imply that $\C_1$ is in class 
$C_{\cdot 0}$ by Proposition~\ref{regC.0}, a contradiction.  There is another
conclusion one can draw from this example, namely that the third
condition in the definition of quasi-free does not follow
automatically from the first two.  In particular, if one considers the
function $1$ as a generator for $H^2(\D)$ and the scalar 1 as a
generator for $\C_1$, the $\mathcal{R}_\mathcal{M}$ space in this case is 
$H^2(D) \oplus \C_1$ and the conditions (i) and (ii) are satisfied. However, the
nonzero vector $f = 1 \oplus 0$ is in $H^2(\D) \oplus \C_1$ but $f
\otimes_\mathcal{A} 1_{z} = 0$ for $z\in \D$.

Another question is whether $\ker Q = \{0\}$ is necessary for
$\mathcal{R}_\mathcal{M}$ to be quasi-free.  However, to establish 
that one would need to
exhibit a nonzero vector $h$ in $\mathcal{R}_\mathcal{M}$ satisfying
$h \otimes_\mathcal{A} 1_{\ul{z}} = 0$ for every $\ul{z}$ in $\Omega$
assuming $\ker Q \not = \{0\}$.  If $h$ is a nonzero vector in $\ker
Q$, then $ (Q \otimes_\mathcal{A} 1_{\ul{z}})(h \otimes_\mathcal{A}
1_{\ul{z}}) = (Qh) \otimes_\mathcal{A} 1_{\ul{z}} = 0$.  The proof
would be completed by showing that $\ker (Q \otimes_{\mathcal{A}}1_{\ul{z}}) =
\{0\}$ for each $\ul{z} \in \Omega$, where $Q \otimes_\mathcal{A}
1_{\ul{z}}: \mathcal{R}_\mathcal{M} \otimes_\mathcal{A} \C_{\ul{z}}
\rightarrow \mathcal{R} \otimes \C_{\ul{z}}$.  The module map $Q
\otimes_\mathcal{A} 1_{\ul{z}}$ is defined by taking the generating
set $\{ (f_i \oplus g_i)\otimes_\mathcal{A} 1_{\ul{z}} \}$ termwise to
the generating set $\{g_i \otimes_\mathcal{A} 1_{\ul{z}} \}$.  If $k <
\infty$, then both of the spaces $\mathcal{R}_\mathcal{M}
\otimes_\mathcal{A} \C_{\ul{z}}$ and $\mathcal{R}
\otimes_\mathcal{A} \C_{\ul{z}}$ are $k$ - dimensional and the map
$Q \otimes_\mathcal{A} 1_{\ul{z}}$ is onto.  Therefore, the kernel 
is trivial and the converse is seen to hold.  For $k = \infty$, we
are unable to conclude that the maps
$Q\otimes_{\mathcal{A}}1_{\ul1{z}}$
have trivial kernels.  To proceed further in the $k = \infty$ case, 
one would need more
information on the nature of the bases defined by the $\{f_i\}$ and
the $\{g_i\}$ and their relationship to each other.

One would like to show that $\ker Q = \{0\}$ if $\mathcal{M}$ belongs to 
class $C_{\cdot 0}$ and a proof would seem almost at hand.  However, 
what the argument seems to requires is assuming that the module 
$\mathcal{M}$ satisfies a stronger condition than that of membership
in the class $C_{\cdot 0}$. One can complete the proof if in the 
definition of class $C_{\cdot 0}$  
the uniform bound on the supremum norm for the sequence in 
$\mathcal{A}(\Omega)$ is  replaced by a 
uniform bound on the Hilbert module norm in the quasi-free module, but 
that would seem to be asking too much.   Thus it is not clear just what 
is the relationship between the notions of class $C_{\cdot 0}$ and weakly 
regular.

Now we come to our principal result, the existence of resolutions.
\begin{Theorem} \label{qfRes}
Every Hilbert module in class {\em (PS)} possesses a weak 
resolution  by quasi-free Hilbert modules.
\end{Theorem}

\begin{proof}  Let us first consider the finitely generated case.  It
is enough to show for $\mathcal{M}$ in class (PS) that the kernel 
of the module map $Q: \mathcal{R}_\mathcal{M}\rightarrow
\mathcal{M}$ is trivial for $\mathcal{R}$ quasi-free.  If
$\mathcal{L}$ is a finite dimensional Hilbert module supported at a
point $\ul{z}$ in $\Omega$, then module multiplication by $\varphi$ on
$\mathcal{L}$ depends only on the values of $\varphi$ and a fixed
finite set of partial derivatives of $\varphi$ at $\ul{z}$.  Now suppose
$Q(k \oplus 0) = 0$ for some vector $k \oplus 0$ in $\mathcal{R}_{\mathcal{M}}$,
which is the form a vector must have in $\ker Q$.  Then there exists a
sequence of functions $\{\varphi_i^{(n)}\}$ such that $\sum \varphi_i^{(n)}
f_i \oplus \varphi_i^{(n)} g_i \rightarrow k \oplus 0$.  By the
definition of the norm on $\mathcal{M} \oplus \mathcal{R}$, and the
fact that $Q(k \oplus 0) = 0$, it follows that $\sum \varphi_i^{(n)} f_i
\rightarrow k$ in $\mathcal{M}$ and $\sum \varphi_i^{(n)} g_i \rightarrow
0$ in $\mathcal{R}$.  Since $\mathcal{R}$ is quasi-free, it follows
that $\sum \varphi_i^{(n)} g_i \rightarrow 0$ implies $\varphi_i^{(n)}
(\ul{z}) \to 0$ and the same is true for evaluation at $\ul{z}$ of the partial
derivatives of $\varphi_i^{(n)}$, all of which follows by localizing
$\mathcal{R}$ with respect to the modules $\mathcal{Q}^n_{\ul{z}}$.  This, of course,
implies that the image 
$\sum \varphi_i^{(n)} g_i \otimes_\mathcal{A} x \rightarrow 0$ 
for $x \in \mathcal{L}$ since 
$\sum \varphi_i^{(n)} g_i \otimes_\mathcal{A} x = 
\sum g_i \otimes_\mathcal{A} \varphi_i^{(n)}(\ul{z})x$ 
and the module action of $\varphi_i^{(n)} (\ul{z})$ on the
vector $x$ in $\mathcal{L}$ depends only on a fixed number of 
partial derivatives of
the functions at $\ul{z}$.  But this implies that the image of $k$ 
is zero in $\mathcal{R} \otimes_{\mathcal{A}}\mathcal{L}$ and hence
by assumption, $k \oplus 0$ is zero or $\ker Q = \{0\}$.

Now suppose $\mathcal{M}$ has infinite rank with generators
$\{f_i\}_{i\in \mathbb{N}}$ and
let $\mathcal{R}$ be a rank one quasi-free Hilbert module over
$\mathcal{A}(\Omega)$ with generator $g$.  Let $\mathcal{M}_i$ be the
submodule of $\mathcal{M}$ generated by $f_i$ for $i \geq 1$.  We can
construct a module $\mathcal{R}_{\mathcal{M}_i}$ for each $i \geq 1$ with
contractive module map $X_i: \mathcal{R}_{\mathcal{M}_i} \rightarrow
\mathcal{M}$ having range dense in $\mathcal{M}_i \subseteq
\mathcal{M}$.  
Since each $\mathcal{M}_i$ is a submodule of $\mathcal{M}$, it follows that
$\mathcal{M}_i$ belongs to class (PS) and hence each
$\mathcal{R}_{\mathcal{M}_i}$ is quasi-free over $\mathcal{A}(\Omega)$
for the basis $\{g \oplus f_i\}$.  Moreover, since the bounds for 
evaluation at $\ul{z}$ on all $\mathcal{R}_{\mathcal{M}_i}$ 
are dominated by the bound on evaluation at 
$\ul{z}$ on $\mathcal{R}$, this implies that
$\mathcal{R}_{\tilde{\mathcal{M}}} = 
\oplus \mathcal{R}_{\mathcal{M}_i}$ is quasi
free over $\mathcal{A}(\Omega)$ for the basis $\{g \oplus f_i:~i \geq 1\}$.
If we define $X: \mathcal{R}_{\tilde{\mathcal{M}}}\rightarrow \mathcal{M}$ such
that $X(\sum k_i) = \sum \frac{1}{2^i} X_i k_i$, then $X$ is 
well-defined and bounded
since $\|X(\sum k_i)\| = \| \sum \frac{1}{2^i} X_i k_i\| \leq \sum \frac{1}{2^i}
\|X_ik_i\| \leq \sum \frac{1}{2^i} \|k_i\| \leq \|\sum k_i\|$.  
This completes the proof of weak regularity in the case of infinite rank.

Given $X: \mathcal{R}_1 \rightarrow \mathcal{M}$  
with $\mathcal{R}_1 = \mathcal{R}_{\mathcal{M}}$ or
$\mathcal{R}_{\tilde{\mathcal{M}}}$, the kernel of 
$X_0=X$ is a submodule of a quasi-free Hilbert module and 
hence belongs to class (PS).  Thus we can repeat the construction 
using  $\ker X_0$ in place of $\mathcal{M}$ and continue
to obtain $X_1:\mathcal{R}_2 \to \mathcal{R}_1$.  
We continue the process using  $\ker X_i$ with the 
process terminating if  $\ker X_i$ is a quasi-free Hilbert module.  
In that case one takes the last $\mathcal{R}_i$ to be the zero module.  
Otherwise, one  continues the process indefinitely obtaining an infinite 
resolution.
\end{proof}
Note that if $\mathcal{M}$ is finitely generated, then the module 
$\mathcal{R}_0$ can be taken to be finitely generated as well. 
However, unless one can conclude that the kernels of the $X_i$ are
finitely generated, we can say nothing about the existence of a
resolution by finitely generated, quasi-free Hilbert modules.
The purpose in proving this result for the class (PS) was to obtain this
finiteness result.  We could have proved the following result
directly which includes Theorem 1.

\begin{Theorem} \label{}
Every Hilbert module in class {\em (SM)} has a weak resolution by
quasi-free Hilbert modules.
\end{Theorem}

\begin{proof}The argument for the infinitely generated case given in the
proof of the preceding theorem can be used, once we know that a singly
generated Hilbert module in class (SM) is weakly regular.  Hence,
assume that $\mathcal{M}$ is a singly generated Hilbert module over 
$\mathcal{A}(\Omega)$ in class (SM) with generating vector $f$, and that 
$\mathcal{R}$ is a singly generated
Hilbert module over $\mathcal{A}(\Omega)$ with generating vector $g$.  
We need to show that
$\ker Q =\{0\}$, where $Q: \mathcal{R}_{\mathcal{M}}\to \mathcal{R}$.  
A vector in the kernel of $Q$ must
have the form $k\oplus 0$ in $\mathcal{M} \oplus \mathcal{R}$.  
Moreover, there exists a sequence
${\varphi_n} \in \mathcal{A}(\Omega)$ such that ${\varphi_n f \oplus \varphi_n g}$ 
converges to $k \oplus 0$,
and hence $\varphi_n f$ converges to $k$ and $\varphi_n g$ converges to $0$.

Now consider the vector $Sk = \lim S(\varphi_n f) = 
\lim (\varphi_n f \otimes_{\mathcal{A}(\Omega)} g) = 
\lim (f \otimes_{\mathcal{A}(\Omega)} \varphi_n g) = 
f \otimes_{\mathcal{A}(\Omega)} (\lim  \varphi_n g) = 0.
$

Therefore, the assumption that $\mathcal{M}$ is smooth for $\mathcal{R}$
and $g$ implies $k = 0$ or $\ker Q= \{0\}$ which completes the proof.
\end{proof}
We can also show that almost compactly supported Hilbert modules have
weak resolutions.

\begin{Theorem} \label{wresqf} 
Every almost compactly supported Hilbert module over
$\mathcal{A}(\Omega)$ has a weak resolution by quasi-free Hilbert modules.
\end{Theorem}

\begin{proof} As before, it is sufficient to show that an almost compactly supported
Hilbert module is weakly regular.  Let $\{f_i\}_{i\in \mathbb{N}}$ be a set of compactly
supported vectors in $\mathcal{M}$ that spans $\mathcal{M}$, and let
$\mathcal{M}_i$ be the submodule of $\mathcal{M}$ generated by $f_i$.
Let $\mathcal{R}$ be the Bergman module
for $\Omega$ with the function $1$ as a basis and let $\mathcal{R}_{\mathcal{M}_i}$
be the module constructed from $\mathcal{R}$ and $\mathcal{M}_i$.
Further, let $X_i$ be the map from $\mathcal{R}_{\mathcal{M}_i}$ to
$\mathcal{M}$ with range dense in $\mathcal{M}_i$.  If we can show that
each $\mathcal{R}_{\mathcal{M}_i}$ is quasi-free, then we can complete
the proof as we did for theorem 1.  Fix an $i \geq 1$.  If $k \oplus 0$ is
in the kernels of the corresponding $Q_i$ from
$\mathcal{R}_{\mathcal{M}_i}$ to $\mathcal{R}$, then there is a
sequence of functions ${\varphi_n} \in \mathcal{A}(\Omega)$ such that
$\varphi_n f_i \oplus \varphi_n 1$ converges to $k \oplus 0$.  In the
Bergman space, it follows that this implies the sequence 
$\{\varphi_n\}_{n\in\mathbb{N}}$
converges uniformly to zero on compact subsets of $\Omega$.  But since the
vector $f_i$ is compactly supported, there exists a compact subset $X$ of
$\Omega$ and a constant $\beta$ such that 
$\|\varphi_n f_i\|_{\mathcal{M}} \leq \beta\|\varphi_n\|_X \|f_i\|_{\mathcal{M}}$ 
and hence $k = \lim \varphi_n f_i = 0$.  Thus $\ker Q_i = \{0\}$ which
completes the proof.
\end{proof}
The purpose of this paper is to establish the existence of weak
resolutions under hypotheses as general as possible.  Unfortunately,
the present theorems are not completely satisfactory in that they do
not provide an intrinsic characterization of those Hilbert modules for
which weak resolutions exist.  In discussing this matter further, let
us focus on the question of when a Hilbert module is weakly regular.
As the proofs indicate, weak regularity rests on the module being
supported in some sense on the open set $\Omega$.  Further consideration of
the module obtained as the quotient of the Hardy module over 
$\mathcal{A}(\D)$ by a submodule determined by a singular inner 
function shows that weak
regularity does not imply almost compactly supported.  Also, this
example shows that while almost compactly supported implies class
$C_{\cdot 0}$, the converse is false.  Thus weak regularity lies
somewhere between almost compactly supported and class $C_{\cdot 0}$, and
perhaps coincides with the latter.

Another question is to determine the class of  modules for which exact
resolutions exist.  The construction presented above seems unlikely to
yield resolutions since that would mean showing that the module map $P:
\mathcal{R}_\mathcal{M} \to \mathcal{M}$ is onto.  Clearly that
depends on having greater control on comparisons between the norms of
the vectors of the form $\sum \varphi_i f_i$ and $\sum \varphi_i g_i$.
Although one could take $\mathcal{R}$ to be the Bergman space, as we
have at various junctures above, the inequalities one would need are
not available, in general.  A different construction, based on the one
given in (Chapter 5, \cite{rgdvip}) might be used to show that a compactly
supported $\mathcal{M}$ is regular but the details are not all clear.
Recall that we do know that a regular Hilbert module lies in class
$C_{\cdot 0}$ but, unfortunately, not the converse.  Understanding
whether class $C_{\cdot 0}$ implies that a Hilbert module is weakly
regular or regular are extremely important in continuing this
approach.

Given that we know so little about regularity, it would seem almost
presumptuous to even introduce the notion of strong regularity and a
modicum of experience would suggest that it hardly ever happens.
However, the resolution provided by the Sz.-Nagy-Foias model is
strongly regular.  Moreover, if one were to consider this question
without appealing to the model theory, we believe it might be hard to
make the right guess.  For example, it would seem unlikely for there
to be a strong resolution of the Bergman module over the disk by a
direct sum of Hardy modules, but there is, albeit one of infinite
multiplicity.  In a different direction, consider the second
resolution of $\C_{\ul{0}}$ given in section 2.  While the maps
are onto, they are not partial isometries.  We presented this
resolution in the form given there because that is the most natural
way.  However, with minor changes in the norms on the resolving
modules, one obtains a strong resolution.  Here are the details.

Recall one maps the Hardy module $H^2(\D^2)$ over the bi-disk algebra
$\mathcal{A}(\D^2)$ onto $\C_{\ul{0}}$ which is a partially isometric
map.  Then one maps the direct sum of two copies of $H^2(\D^2)$ to
$H^2(\D^2)$ by the map $X (f \oplus g) = z_2f - z_1g$ which is not a
partial isometry.  However, suppose one changes the norms on the Hardy
modules so that on the first one, if $a_{i,j}$ are the Taylor
coefficients of $f$, then we multiply the $a_{i,0}$ by a factor
$\scriptstyle{\sqrt{1/2}}$ and on the second one, if $b_{i,j}$ are the
Taylor coefficients of $g$, then we multiply the $b_{0,j}$ by a factor
$\scriptstyle{\sqrt{1/2}}$.  The resulting Hilbert modules are still
quasi-free over $\mathcal{A}(\D^2)$ since both changes yield 
equivalent norms.  However, now $X$ is a partial isometry.  Now, the
last non-zero term in the resolution is a copy of $H^2(\D^2)$ with the
map $Y$ defined by $Yf = z_1f \oplus z_2f$.  Here,
one wants to multiply all the Taylor coefficients of a vector $f$ in
$H^2(\D^2)$ by a factor $\scriptstyle{\sqrt{1/2}}$ to obtain an
equivalent Hilbert module which is quasi-free and which makes $Y$ into
a partial isometry.  Thus, $\C_{\ul{0}}$ has a strong resolution.

The question of whether or not resolutions or strong resolutions exist
is not merely academic for the following reason.  If a Hilbert module
$\mathcal{M}$ is regular, then there exists a quasi-free Hilbert
module $\mathcal{R}$ and a module map $X$ from $\mathcal{R}$ onto
$\mathcal{M}$.  If $\ker X = \{0\}$, then $\mathcal{M}$ is similar to the
quasi-free module $\mathcal{R}$ and hence, $\mathcal{R}$ is itself quasi-free.
Otherwise, we may assume there is a nontrivial kernel.  If a full
resolution exists, then there are nontrivial kernels and we can
continue or the resolution stops and has finite length.  This is the
situation in commutative algebra and one should expect in such a case
to be able to extract information about $\mathcal{M}$ from the
resolution using the extension of techniques from commutative algebra.

There is another issue which it is convenient to raise at this time,
namely, are resolutions finite?  In general, the answer must be negative.
However, one would expect that there is a large class of Hilbert
modules for which that is the case.  A related question is whether the
kernel of a module map $X: \mathcal{R} \to \mathcal{M}$ is
finitely generated.  Again, one would assume that this is the case for
a large class of Hilbert modules when both $\mathcal{M}$ and $\mathcal{R}$ are
finitely generated, but results seem to be difficult \cite{yang}.
The questions we are raising here, of course,  
concern coherence-like properties of the spectral sheaf 
${\rm Sp}(\mathcal{M})$.  If one replaces
Hilbert modules by Frechet modules, then there is a lot of work on
these questions (cf. \cite{EM}).  Connecting the two approaches, where
Hilbert spaces are used on the one hand or Frechet spaces on the other,
seems difficult.  Our point of view has been that the appropriate
domain for multivariate operator theory is Hilbert space but any final
assessment must rest on the results obtained and their utility.

\section{Usefulness of Resolutions}

This paper has been devoted to showing the existence of quasi-free
resolutions of Hilbert modules.  There would be little point in
constructing such resolutions if they were not useful in studying the
original modules.  In this section we want to sketch some ways in
which resolutions have been useful and could be useful in the general
study of Hilbert modules.

As we have indicated, one can re-interpret the model of Sz.-Nagy and
Foias as a resolution.  Thus, in principle, one could argue that all
of model theory could be taken over to the context of resolutions but
that would be an exaggeration.  Much of the theory depends on the rich
interplay of function theory, functional analysis, and Fourier
analysis which come together on the unit disk with boundary the unit
circle.  Also, some of the theory depends on the fact one has a strong
resolution in our language rather than just a resolution or even a
weak resolution.  Still the characterization of at least one basic
notion carries over, that of the spectrum.

Recall that the resolution for a contractive Hilbert module H of class 
$C_{\cdot 0}$ has the form:
$$
0 \longleftarrow  \mathcal{H} \longleftarrow  H^2_{\mathcal{E}}(\D)  
\stackrel{X}{\longleftarrow}H^2_{\mathcal{E}^*}(\D) \longleftarrow 0,
$$ 
where $X$ is an isometric module map.  If one localizes $X$ by $C_z$, 
one can show that $X \otimes_{\mathcal{A}(\mathbb{D})} 1_z = 
\Theta(z)$, where $\Theta$ is the characteristic
operator function of Sz.-Nagy and Foias. In general model theory, 
one knows that $\Theta$ is an
operator-valued inner function, that is, it has radial limits a.\ e.\ on
the circle which are unitary operators from $\mathcal{E}^*$ to
$\mathcal{E}$.   Our interest is
in the connection of $\Theta$ with the spectrum which one knows is the
union of the set of points in $\D$ at which $\Theta(z)$ fails to be
invertible plus the closed subset of the boundary on which $\Theta$ fails
to have an analytic continuation.  One can show directly from the
exactness of the resolution that the spectrum inside $\D$ consists of the
points at which the localization $X \otimes_{\mathcal{A}(\mathbb{D})} 1_z$ 
fails to be invertible and, indeed, that the nature of the spectrum 
of the operator defined by module multiplication by $z$ is that 
same as that of $\Theta(z)$.  The details of this calculation are
given in \cite{rgdvip}.

Now suppose we have a weak resolution of a Hilbert module
$\mathcal{M}$ over $\mathcal{A}(\Omega)$.  One can calculate the
spectrum of the module which is defined using the Taylor spectrum
(cf. \cite{rgdvip}), in terms of the resolution.  Moreover, one can
determine the nature of the spectrum, that is, the nature of the 
lack of exactness of the Koszul complex a la Taylor.  
One should compare a recent paper by D.
Greene \cite{DG} in which he does something similar for modules over an
algebra of holomorphic functions but one which is not a function
algebra.  In both cases, the behavior of the Hilbert module on the
boundary would have to be investigated using different techniques.  As we
indicated above, on the disk the determination of the full 
spectrum involves the notion of analytic continuation.
Although, there are other characterizations,  
none involve strictly algebraic notions.

There is another class of invariants for Hilbert modules 
of a very different nature,
associated with complex geometry.  In the late seventies, M. Cowen and
the first author introduced a class of operators which have a
hermitian holomorphic bundle associated with them. Moreover, they showed that the
geometric invariants of the bundle form a complete set of unitary
invariants for the operator.  This approach was extended by Curto and
Salinas \cite{rcns} to the case of commuting n-tuples of operators,
and by X. Chen and the first author \cite{chenrgd} to certain classes
of Hilbert modules.  The latter class includes the quasi-free ones and
the associated spectral sheaves are the corresponding hermitian
holomorphic vector bundles.  Thus quasi-free Hilbert modules can be
characterized up to unitary equivalence by the curvature and a finite
set of partial derivatives of curvature over $\Omega$.  In a series of
papers \cite{rgdgm}, \cite{DM}, \cite{DMC}, \cite{DMCQ}, the authors along with
Verughese, have related the geometrical invariants for Hilbert modules
in a resolution.  In particular, one shows for the quotient module
defined by the functions in a quasi-free module $\mathcal{R}$ that
vanish to some order along a hypersurface, that the geometric
invariants for the spectral sheaf for the quotient are determined by
those for the quasi-free sheaf in the form of longitudinal curvature,
transverse curvature and a second fundamental form involving an
appropriate jet bundle.

One can formulate relations such as the above for weak resolutions
although the formulas and proofs will involve, ultimately, an
extension of techniques related to the work of Harvey-Lawson \cite{hl}
as well as to that of Demailley \cite{De}.  Some very simple cases
have been established but there is much to do and the possibility for
relating unitary invariants for a module to those of a weak resolution
seem promising.

\end{document}